\documentclass[10pt]{article}
\begin{document}
\title{ {\bf Potentially $K_{m}-G$-graphical Sequences: A Survey }
\thanks{  Project Supported by  NSF of Fujian(2008J0209),
 Fujian Provincial
Training Foundation for "Bai-Quan-Wan Talents Engineering" , Project
of Fujian Education Department and Project of Zhangzhou Teachers
College (SK07009).}}
\author{{Chunhui Lai, Lili Hu}\\
{\small Department of Mathematics and Information Science,} \\
{\small Zhangzhou Teachers College,}
\\{\small Zhangzhou, Fujian 363000,
 P. R. of CHINA.}\\{\small  zjlaichu@public.zzptt.fj.cn (Chunhui
 Lai, Corresponding author)}
 \\{\small  jackey2591924@163.com ( Lili Hu)}}
\date{}
\maketitle
\begin{center}
\begin{minipage}{4.2in}
\vskip 0.1in
\begin{center}{\bf Abstract}\end{center}
 {The set of all non-increasing nonnegative integers sequence $\pi=$
  ($d(v_1 ),$ $d(v_2 ),$ $...,$ $d(v_n )$) is denoted by $NS_n$.
  A sequence
$\pi\in NS_n$ is said to be graphic if it is the degree sequence of
a simple graph $G$ on $n$ vertices, and such a graph $G$ is called a
realization of $\pi$. The set of all graphic sequences in $NS_n$ is
denoted by $GS_n$. A graphical sequence $\pi$ is potentially
$H$-graphical if there is a realization of $\pi$ containing $H$ as a
subgraph, while $\pi$ is forcibly $H$-graphical if every realization
of $\pi$ contains $H$ as a subgraph. Let $K_k$ denote a complete
graph on $k$ vertices. Let $K_{m}-H$ be the graph obtained from
$K_{m}$ by removing the edges set $E(H)$ of the graph $H$ ($H$ is a
subgraph of $K_{m}$). This paper summarizes briefly some recent
results on potentially $K_{m}-G$-graphic sequences and give a useful
classification for determining $\sigma(H,n)$.}
  \par
\end{minipage}
\end{center}
 \par
 \section{Introduction}
\par

  The set of all non-increasing nonnegative integer sequence $\pi=$
  ($d(v_1 ),$ $d(v_2 ),$ $...,$ $d(v_n )$) is denoted by $NS_n$.
  A sequence
$\pi\in NS_n$ is said to be graphic if it is the degree sequence of
a simple graph $G$ on $n$ vertices, and such a graph $G$ is called a
realization of $\pi$. The set of all graphic sequences in $NS_n$ is
denoted by $GS_n$. A graphical sequence $\pi$ is potentially
$H$-graphical if there is a realization of $\pi$ containing $H$ as a
subgraph, while $\pi$ is forcibly $H$-graphical if every realization
of $\pi$ contains $H$ as a subgraph. If $\pi$ has a realization in
which the $r+1$ vertices of  largest degree induce a clique, then
$\pi$ is said to be potentially $A_{r+1}$-graphic. Let
$\sigma(\pi)=d(v_1 )+d(v_2 )+... +d(v_n ),$ and $[x]$ denote the
largest integer less than or equal to $x$. We denote $G+H$ as the
graph with $V(G+H)=V(G)\bigcup V(H)$ and $E(G+H)=E(G)\bigcup
E(H)\bigcup \{xy: x\in V(G) , y \in V(H) \}. $ Let $K_k$, $C_k$,
$T_k$, and $P_{k}$ denote a complete graph on $k$ vertices,  a cycle
on $k$ vertices, a tree on $k+1$ vertices, and a path on $k+1$
vertices, respectively. Let $F_k$ denote the friendship graph on
$2k+1$ vertices, that is, the graph of $k$ triangles intersecting in
a single vertex. For $0\leq r \leq t,$ denote the generalized
friendship graph on $kt-kr+r$ vertices by $F_{t,r,k}$, where
$F_{t,r,k}$ is the graph of $k$ copies of $K_t$ meeting in a common
$r$ set. We use the symbol $Z_4$ to denote $K_4-P_2.$  Let $K_{m}-H$
be the graph obtained from $K_{m}$ by removing the edges set $E(H)$
of the graph $H$ ($H$ is a subgraph of $K_{m}$).
\par

Given a graph $H$, what is the maximum number of edges of a graph
with $n$ vertices not containing $H$ as a subgraph? This number is
denoted $ex(n,H)$, and is known as the Tur\'{a}n number. In terms of
graphic sequences, the number $2ex(n,H)+2$ is the minimum even
integer $l$ such that every $n$-term graphical sequence $\pi$ with
$\sigma (\pi) \geq l $ is forcibly $H$-graphical. Erd\"os,\ Jacobson
and Lehel [13] first consider the following variant: determine the
minimum even integer $l$ such that every $n$-term graphical sequence
$\pi$ with $\sigma(\pi)\ge l$ is potentially $H$-graphical. We
denote this minimum $l$ by $\sigma(H, n)$. Erd\"os,\ Jacobson and
Lehel [13] showed that $\sigma(K_k, n)\ge (k-2)(2n-k+1)+2$ and
conjectured that equality holds. They proved that if $\pi$ does not
contain zero terms, this conjecture is true for $k=3,\ n\ge 6$. The
conjecture is confirmed in [19] and [43-46]. Li et al. [46] and
Mubayi [55] also independently determined the values
$\sigma(K_r,2k)$ for any $k \geq 3$ and $n\geq k$. Li and Yin [51]
further determined $\sigma(K_r,n)$ for $r\geq7$ and $n\geq2r+1$. The
problem of determining $\sigma(K_r,n)$ is completely solved.

\par
  Gould, Jacobson and Lehel[19] also proved that
  $\sigma(pK_2,n)=(p-1)(2n-p)+2$ for $p\geq2$; $\sigma(C_4,n)=2[{{3n-1}\over 2}]$ for $n\ge 4$. Lai[29]
     gave a lower bound of $\sigma(C_k,n)$ and proved that
     $\sigma(C_5,n)=4n-4$ for $n\geq5$ and $\sigma(C_6,n)=4n-2$ for
     $n\geq7$. Lai [32] proved that $\sigma(C_{2m+1}, n)=m(2n-m-1)+2$,
for $m\geq 2, n\geq 3m$; $\sigma(C_{2m+2} , n)=m(2n-m-1)+4$,
 for $ m\geq 2, n\geq 5m-2$. Li and Luo[41] gave a lower bound of $\sigma(_3C_l,n)$
 and determined $\sigma(_3C_l,n)$, $4\leq l\leq6, n\geq l$. Li, Luo and
 Liu [42] determined $\sigma(_3C_l,n)$ for $ 3\leq l\leq 8$ and $n\geq l.$ and
 $\sigma(_3C_9,n)$ for  $n\geq 12$. Li and Yin [48] determined $\sigma(_3C_l,n)$
 for $n$ sufficiently large.  Yin, Li
 and Chen[68] determined $\sigma(_kC_l,n)$, $l\geq7, 3\leq k\leq l$.
 Chen and Yin[9] determined the
 values $\sigma(W_5,n)$ for $n\geq11$ where $W_r$ is a wheel graph
 on $r$ vertices.  For $r\times s$ complete bipartite graph $K_{r,s}$,
Gould, Jacobson and Lehel[19] determined $\sigma(K_{2,2},n)$.   Yin
et al. [63,65,69,70] determined $\sigma(K_{r,s},n)$ for $s\geq r\geq
2$ and sufficiently large $n$. For $r\times s \times t$ complete
3-partite graph $K_{r,s,t}$, Erd\"{o}s, Jacobson and Lehel[13]
determined $\sigma(K_{1,1,1},n)$. Lai[30] determined
$\sigma(K_{1,1,2},n)$. Yin[58] and Lai[34] independently determined
$\sigma(K_{1,1,3},n)$.
       Chen[7] determined $\sigma(K_{1,1,t},n)$ for $t\geq3$, $n\geq2[{{(t+5)^2}\over
     4}]+3$. Chen[5] determined $\sigma(K_{1,2,2},n)$ for $5\leq
     n\leq8$ and $\sigma(K_{2,2,2},n)$ for $n\geq6$.
      Let $K_s^t$ denote the complete $t$ partite graph such that each partite
  set has exactly $s$ vertices.
  Guantao Chen, Michael Ferrara, Ronald J.Gould, John R. Schmitt[11] showed that
  $\sigma(K_s^t,n)=\pi(K_{(t-2)s}+K_{s,s},n)$ and obtained  the exact value of $\sigma(K_j+K_{s,s},n)$
 for $n$ sufficiently large. Consequently, they obtained the exact value of $\sigma(K_s^t,n)$ for $n$ sufficiently large.
  For $n\geq5$,    Ferrara, Jacobson and Schmitt[17] determined $\sigma(F_k,n)$
where $F_k$ denotes the graph of $k$ triangles intersecting at
exactly one common vertex. In[16], Ferrara, Gould and Schmitt
determined a lower bound for $\sigma(K_s^t,n)$ where $K_s^t$ denotes
the complete multipartite graph with $t$ partite sets each of size
$s$, and proved equality in the case $s=2$. They also provided a
graph theoretic proof of the value of $\sigma(K_t,n)$. Michael J.
Ferrara[15] determined $\sigma(H,n)$ for the graph $H=K_{m_1}\cup
K_{m_2}\cup \cdots\cup K_{m_k}$ where $n$ is sufficiently large
integer. Ferrara, M., Jacobson, M., Schmitt, J. and Siggers M.[18]
determined $\sigma(K_{s,t},m,n)$, $\sigma(P_t,m,n)$ and
$\sigma(C_{2t},m,n)$ where $\sigma(H,m,n)$ is the minimum integer
$k$ such that every bigraphic pair $S=(A,B)$ with $|A|=m$, $|B|=n$
and $\sigma(S)\geq k$ is potentially $H$-bigraphic. For an
arbitrarily chosen $H$, Schmitt, J.R. and Ferrara, M.[56] gave a
good lower bound of $\sigma(H,n)$.   Yin and Li[67] determined
$\sigma(K_{r_1,r_2,\cdots,r_l,r,s,}n)$ for sufficiently large $n$.
Moreover, Yin, Chen and Schmitt[62] determined $\sigma(F_{t,r,k},n)$
for $k\geq2,t\geq3,1\leq r\leq t-2$ and sufficiently large, where
$F_{t,r,k}$ denotes the graph of $k$ copies of $K_t$ meeting in a
common $r$ set.  Gupta, Joshi and Tripathi [20] gave a necessary and
sufficient condition for the existence of a tree of order
 $n$ with a given degree set. Yin [59] gave a new necessary and sufficient condition for
 $\pi$ to be potentially $K_{r+1}$- graphic. Jiong-sheng Li and Jianhua Yin [50] gave a survey
on graphical sequences.
  \par
  \section{Potentially $K_{m}-G$-graphical Sequences}
   \par
Let $H$ be a graph with $m$ vertices, then $ H = K_m - (K_m - H ).$
Let $G = K_m - H,$ then $ \sigma (H, n) = \sigma (K_m-G, n).$ If the
Problem 1 - 5 in the Open Problems section be solved, then the
problem of determining $\sigma (H, n)$ is completely solved. We
think the Problem 3, 4 is a useful classification for determining
$\sigma(H, n)$.
   \par
     Gould, Jacobson and Lehel[19]
 pointed out that it would be nice to see where in the range for
$3n-2$  to $4n-4,$ the value  $\sigma (K_4-e, n)$ lies. Later,
Lai[30] proved that
\par
{\bf  Theorem 1.} For $n=4,5$ and $n\geq 7$
 \[\sigma(K_4-e,n)=\left\{ \begin{array}{cc}3n-1 & \mbox {if $n$
is odd}  \\3n-2 & \mbox {if $n$ is even.}\end{array}\right.\] For $
n=6,$ $S$ is a 6-term graphical sequence with $\sigma(S) \geq 16$,
then either there is a  realization of S containing $K_4-e$ or
$S=(3^6)$. (Thus $\sigma(K_4-e, 6)=20)$.
\par
      Huang[26] gave a lower bound of $\sigma(K_m-e,n)$.
      Yin, Li and Mao[71] and Huang[27] independently determined the
 values $\sigma(K_5-e,n)$ as following.
 \par
 {\bf  Theorem 2.} If  $n\geq 5$, then $$ \sigma (K_5-e, n) =\left\{
    \begin{array}{ll}5n-7, \ \ \  \mbox{ if $n$ is odd}\\
    5n-6,
     \ \ \ \  \mbox{if $n$ is even} \end{array} \right. $$
 Lai[35-36] determined $\sigma (K_{5}-C_{4}, n),\sigma (K_{5}-P_{3}, n)$ and
    $\sigma (K_{5}-P_{4}, n)$.
    \par
    {\bf  Theorem 3.} For  $n\geq 5$,
    $\sigma (K_{5}-C_{4}, n)=\sigma (K_{5}-P_{3}, n)=\sigma (K_{5}-P_{4}, n)=4n-4.$
    \par
    Yin and Li[66] gave a good method (Yin-Li method) of
determining the values $\sigma(K_{r+1}-e,n)$ (In fact, Yin and
Li[66] also determining the values $\sigma(K_{r+1}-ke,n)$ for
$r\geq2$ and $n\geq 3r^2-r-1$).
\par
{\bf Theorem 4.}
     Let $n\geq r+1$ and $\pi=(d_1,d_2,\cdots,d_n)\in
    GS_n$ with $d_{r+1}\geq r$. If $d_i\geq 2r-i$ for
    $i=1,2,\cdots,r-1$, then $\pi$ is potentially $A_{r+1}$-graphic.

{\bf Theorem 5.}
     Let $n\geq 2r+2$ and $\pi=(d_1,d_2,\cdots,d_n)\in
    GS_n$ with $d_{r+1}\geq r$. If $d_{2r+2}\geq r-1$ , then $\pi$ is
    potentially $A_{r+1}$-graphic.

{\bf Theorem 6.}
     Let $n\geq r+1$ and $\pi=(d_1,d_2,\cdots,d_n)\in
    GS_n$ with $d_{r+1}\geq r-1$. If $d_i\geq 2r-i$ for
    $i=1,2,\cdots,r-1$, then $\pi$ is potentially $K_{r+1}-e$-graphic.

{\bf Theorem 7.}
     Let $n\geq 2r+2$ and $\pi=(d_1,d_2,\cdots,d_n)\in
    GS_n$ with $d_{r-1}\geq r$. If $d_{2r+2}\geq r-1$ , then $\pi$ is
    potentially $K_{r+1}-e$
    -graphic.

{\bf  Theorem 8.}
 If $r \geq 2$ and $n\geq 3r^2-r-1,$ then
  $$\sigma (K_{r+1}-ke, n) =\left\{
    \begin{array}{ll}(r-1)(2n-r)-(n-r)+1, \\ \mbox{ if $n-r$ is odd}\\
    (r-1)(2n-r)-(n-r)+2,
     \\ \mbox{if $n-r$ is even} \end{array} \right. $$
     \par
    After reading[66], using Yin-Li method
    Yin[72] determined the values $\sigma(K_{r+1}-K_3,n)$ for
    $r\geq3,n\geq3r+5$.
    \par
    {\bf  Theorem 9.}
 If $r \geq 3$ and $n\geq 3r+5,$ then
  $\sigma (K_{r+1}-K_3, n) =
    (r-1)(2n-r)-2(n-r)+2.$
      \par
     Determining $\sigma(K_{r+1}-H,n)$, where $H$
    is a tree on 4 vertices is more useful than a cycle on 4
    vertices (for example, $C_4 \not\subset C_i$, but $P_3 \subset C_i$ for $i\geq 5$).
    So, after reading[66] and [72], using Yin-Li method Lai and Hu[38] determined
    $\sigma(K_{r+1}-H,n)$ for  $n\geq 4r+10, r\geq 3, r+1 \geq k \geq 4$ and $H$ be a graph on $k$
    vertices which
    containing a tree on $4$ vertices but
     not containing a cycle on $3$ vertices and $\sigma (K_{r+1}-P_2, n)$ for
    $n\geq 4r+8, r\geq 3$.
    \par
    {\bf  Theorem 10.} If $r\geq3$ and $n\geq 4r+8$, then
          $\sigma(K_{r+1}-P_2,n)=(r-1)(2n-r)-2(n-r)+2.$

\par
{\bf  Theorem 11.} If $r\geq3, r+1 \geq k \geq 4$ and $n\geq 4r+10$,
then $\sigma (K_{r+1}-H, n)= (r-1)(2n-r)-2(n-r),$ where $H$ is a
graph on $k$
    vertices which contains a tree on $4$ vertices but
     not contains a cycle on $3$ vertices.
       \par

     There are a number of graphs on $k$
    vertices which
    containing a tree on $4$ vertices but
     not containing a cycle on $3$ vertices (for example, the cycle
      on $k$ vertices, the tree on $k$ vertices, and the
   complete 2-partite graph on $k$ vertices, etc ).
   \par
   Using Yin-Li method Lai and Sun[39] determined
$\sigma (K_{r+1}-(kP_2\bigcup tK_2), n)$ for
    $n\geq 4r+10, r+1 \geq 3k+2t,
    k+t \geq 2,k \geq 1, t \geq 0$.
    \par
{\bf  Theorem 12.} If $n\geq 4r+10,  r+1 \geq 3k+2t,
    k+t \geq 2,k \geq 1, t \geq 0$, then $\sigma (K_{r+1}-(kP_2\bigcup tK_2), n)=
(r-1)(2n-r)-2(n-r)$.
 \par
 To now, the problem of determining $\sigma
(K_{r+1}-H, n)$ for $H$ not containing a cycle on 3 vertices and
sufficiently large $n$ has been solved.
\par
 Using Yin-Li method Lai[37] determined
    $\sigma(K_{r+1}-Z,n)$ for $n\geq5r+19, r+1\geq k\geq5, j\geq5$
    and $Z$ is a graph on $k$ vertices and $j$ edges which contains
    a graph $Z_4$ but not contain a cycle on 4 vertices. In the same
    paper, the author also determined the values of  $\sigma (K_{r+1}-Z_4, n)$, $\sigma (K_{r+1}-(K_4-e), n)$,
      $\sigma (K_{r+1}-K_4, n)$ for
    $n\geq 5r+16, r\geq 4$.
    \par
    {\bf  Theorem 13.} If $r\geq 4$ and $n\geq 5r+16$, then
 $$ \sigma (K_{r+1}-K_{4}, n) = \sigma (K_{r+1}-(K_{4}-e), n) =$$
 $$\sigma (K_{r+1}-Z_{4}, n) =\left\{
    \begin{array}{ll}(r-1)(2n-r)-3(n-r)+1, \\ \mbox{ if $n-r$ is odd}\\
    (r-1)(2n-r)-3(n-r)+2,
     \\ \mbox{if $n-r$ is even} \end{array} \right. $$

\par
{\bf  Theorem 14.} If $n\geq 5r+19,  r+1 \geq k \geq 5,$ and $j \geq
5$, then
  $$\sigma (K_{r+1}-Z, n) =\left\{
    \begin{array}{ll}(r-1)(2n-r)-3(n-r)-1, \\ \mbox{ if $n-r$ is odd}\\
    (r-1)(2n-r)-3(n-r)-2,
     \\ \mbox{if $n-r$ is even} \end{array} \right. $$
      where $Z$ is a graph on $k$
    vertices and $j$ edges which
    contains a graph $Z_4$ but
     not contains a cycle on $4$ vertices.
    \par

There are a number of graphs on $k$
    vertices  and $j$ edges which
    contains a graph $Z_4$  but
     not contains a cycle on $4$ vertices.
     (for example,  the graph obtained by $C_3, C_{i_1}, C_{i_2},
     \cdots, C_{i_p}$ intersecting in a single vertex
     $(i_j\neq 4, j=1,2,3,\cdots, p)$ (if $i_j=3, j=1,2,3,\cdots, p,$
     then this graph is the friendship graph $F_{p+1}$, ),
     the graph obtained by
     $C_3, P_{i_1}, P_{i_2}, \cdots, P_{i_p}$ intersecting in a single
     vertex,
     $(i_1 \geq 1)$, the graph obtained by $C_3, P_{i_1}, C_{i_2},
     \cdots, C_{i_p}$ $(i_j\neq 4, j=2,3,\cdots, p, i_1 \geq 1)$
     intersecting in a single vertex, etc )
       \par
       \par
    Using Yin-Li method Lai and Yan[40] proved that
    \par
    {\bf  Theorem 15.}
 If $n\geq 5r+18,  r+1 \geq k \geq 7,$ and $j
\geq 6$, then
  $$\sigma (K_{r+1}-U, n) =\left\{
    \begin{array}{ll}(r-1)(2n-r)-3(n-r)-1, \\ \mbox{ if $n-r$ is odd}\\
    (r-1)(2n-r)-3(n-r),
     \\ \mbox{if $n-r$ is even} \end{array} \right. $$
      where $U$ is a graph on $k$
    vertices and $j$ edges which
    contains a graph $(K_{3} \bigcup P_{3})$ but
     not contains a cycle on $4$ vertices and not contains $Z_4$.

There are a number of graphs on $k$
    vertices  and $j$ edges which
    contains a graph $(K_{3} \bigcup P_{3})$  but
     not contains a cycle on $4$ vertices and not contains $Z_4$.
     (for example, $C_3\bigcup C_{i_1} \bigcup C_{i_2} \bigcup
     \cdots \bigcup C_{i_p}$ $(i_j\neq 4, j=2,3,\cdots, p, i_1 \geq 5)$,
     $C_3\bigcup P_{i_1} \bigcup P_{i_2} \bigcup \cdots \bigcup P_{i_p}$
     $(i_1 \geq 3)$, $C_3\bigcup P_{i_1} \bigcup C_{i_2} \bigcup
     \cdots \bigcup C_{i_p}$ $(i_j\neq 4, j=2,3,\cdots, p, i_1 \geq 3)$, etc )
       \par

\par
   A harder question is to characterize the potentially
 $H$-graphic sequences without zero terms.  Luo [53] characterized the potentially
 $C_k$-graphic sequences for each $k=3,4,5$.
 \par
 \textbf{\noindent Theorem 16.}  Let $\pi=(d_1,d_2,\cdots,d_n)$
 be a graphic sequence with $n\geq 3$. Then $\pi$ is potentially
$C_3$-graphic if and only if $d_3 \geq 2 $ except for 2 case: $\pi =
(2^4)$ and  $\pi = (2^5)$.

\par
\textbf{\noindent Theorem 17.}  Let $\pi=(d_1,d_2,\cdots,d_n)$
 be a graphic sequence with $n\geq 4$. Then $\pi$ is potentially
$C_4$-graphic if and only if the following conditions hold:
\par
  $(1)$ $d_4\geq 2$.
  \par
  $(2)$ $d_1 =n-1$ implies $d_2\geq 3$
\par
  $(3)$ If $n = 5, 6,$ then $\pi\neq (2^n)$.
\par
\textbf{\noindent Theorem 18.}  Let $\pi=(d_1,d_2,\cdots,d_n)$
 be a graphic sequence with $n\geq 5$. Then $\pi$ is potentially
$C_5$-graphic if and only if the following conditions hold:
\par
  $(1)$ $d_5\geq 2$.
  \par
  $(2)$ For $i = 1, 2, d_1 =n-i$ implies $d_{4-i}\geq 3$
\par
  $(3)$ If $\pi= (d_1, d_2, 2^k, 1^{n-k-2}),$ then $d_1+ d_2 \leq n+k-2$.
\par

Chen [2] characterized the potentially
 $C_k$-graphic sequences for each $k=6$.
 \par
 \textbf{\noindent Theorem 19.}  Let $\pi=(d_1,d_2,\cdots,d_n)$
 be a graphic sequence with $n\geq 6$. Then $\pi$ is potentially
$C_6$-graphic if and only if the following conditions hold:
\par
  $(1)$ $d_6\geq 2$.
  \par
  $(2)$ If $n = 7, 8,$ then $\pi\neq (2^n)$.
  \par
  $(3)$ For $i = 1, 2, 3, d_1 =n-i$ implies $d_{5-i}\geq 3$
\par
  $(4)$ If $\pi= (d_1, d_2, 2^k, 1^{n-k-2}),$ then $d_1+ d_2 \leq n+k-2$;
  if $\pi= (d_1, d_2, 3,  2^k, 1^{n-k-3}),$ then $d_1+ d_2 \leq n+k$;
  if $\pi= (d_1, d_2, 3, 3, 2^k, 1^{n-k-4}),$ then $d_1+ d_2 \leq n+k+2$.
\par
 Yin, Chen and Chen [60] characterized the potentially
 $\sb kC\sb l$-graphic sequences for each $k=3, 4\leq l \leq 5$ and $k=4, l=5$.
 \par
  \textbf{\noindent Theorem 20.}  Let $\pi=(d_1,d_2,\cdots,d_n) \in
  GS_n$
 be a potentially $C_4$-graphic sequence. Then $\pi$ is potentially
$\sb 3 C\sb 4$-graphic if and only if $\pi$ satisfies one of the
following conditions:
\par
  $(1)$ $d_2\geq 3$ and $\pi\neq (3^2, 2^4)$;
  \par
  $(2)$ $\pi= (d_1,  2^k, 1^{n-k-1})$ with $2 \leq d_1 \leq 3$ and $k \geq 6,$
  and  $\pi\neq (2^8)$ and $(2^9)$;
  \par
  $(3)$ $\pi= (d_1,  2^k, 1^{n-k-1})$ with $4 \leq d_1 \leq n-2$ and $k \geq 5,$
  and  $\pi\neq (4, 2^6)$ and $(4, 2^7)$.

\par
 \textbf{\noindent Theorem 21.}  Let $\pi=(d_1,d_2,\cdots,d_n) \in
  GS_n$
 be a potentially $C_5$-graphic sequence. Then $\pi$ is potentially
$\sb 3 C\sb 5$-graphic if and only if $\pi$ satisfies one of the
following conditions:
\par
  $(1)$ $d_2\geq 3$ and $\pi\neq (3^2, 2^4)$ and $(3^2, 2^5)$;
  \par
  $(2)$ $\pi= (d_1,  2^k, 1^{n-k-1})$ with $2 \leq d_1 \leq 3$ and $k \geq 11,$
  and  $\pi\neq (2^{13})$ and $(2^{14})$;
  \par
  $(3)$ $\pi= (d_1,  2^k, 1^{n-k-1})$ with $4 \leq d_1 \leq 5$ and $k \geq 10,$
  and  $\pi\neq (4, 2^{11})$ and $(4, 2^{12})$.
\par
  $(4)$ $\pi= (d_1,  2^k, 1^{n-k-1})$ with $6 \leq d_1 \leq n-4$ and $k \geq 9,$
  and  $\pi\neq (4, 2^{10})$ and $(4, 2^{11})$.
\par
\par
 \textbf{\noindent Theorem 22.}  Let $\pi=(d_1,d_2,\cdots,d_n) \in
  GS_n$
 be a potentially $C_5$-graphic sequence. Then $\pi$ is potentially
$\sb 4 C\sb 5$-graphic if and only if $\pi$ satisfies one of the
following conditions:
\par
  $(1)$ $d_2\geq 3$;
  \par
  $(2)$ $\pi= (d_1,  2^k, 1^{n-k-1})$ with $2 \leq d_1 \leq 3$ and $k \geq 8,$
  and  $\pi\neq (2^{10})$ and $(2^{11})$;
  \par
  $(3)$ $\pi= (d_1,  2^k, 1^{n-k-1})$ with $4 \leq d_1 \leq n-4$ and $k \geq 7,$
  and  $\pi\neq (4, 2^{8})$ and $(4, 2^{9})$.
\par

 Chen, Yin and Fan [10] characterized the potentially
 $\sb kC\sb l$-graphic sequences for each $3 \leq k \leq 5, l=6$.
 \par
 \par
 \textbf{\noindent Theorem 23.}  Let $\pi=(d_1,d_2,\cdots,d_n) \in
  GS_n$, $n \geq 6,$ and $\pi\neq (3^2, 2^{10})$,$(2^{19})$, $(2^{20})$, $(4,
  2^{17})$, $(4,  2^{18})$, $(6,  2^{16})$, $(6,  2^{17})$, $(8,  2^{15})$, $(8,
  2^{16})$,
  Then $\pi$ is potentially
$\sb 3 C\sb 6$-graphic if and only if $\pi$ be a potentially
$C_6$-graphic sequence, and  $\pi$ satisfies one of the following
conditions:
\par
  $(1)$ $d_3\geq 3$, and if $d_1 =d_3 =3, d_4 =2,$ then $d_{10} =2$;
  \par
  $(2)$ $d_2 \geq 4, d_3 =2, d_7 = 2$;
  \par
  $(3)$ $d_2 =3, d_3 =2,$ and if  $4 \geq d_1 \geq 3,$ then $d_{10} =2,$ and if  $n-4 \geq d_1 \geq 5,$ then $d_{9} =2$.
\par
  $(4)$ $d_2 =2,$ and if  $3 \geq d_1 \geq 2,$ then $d_{18} =2,$ and if  $5 \geq d_1 \geq 4,$ then $d_{17} =2$,
  and if  $7 \geq d_1 \geq 6,$ then $d_{16} =2,$ and if  $n-7 \geq d_1 \geq 8,$ then $d_{15} =2$.
\par
 \textbf{\noindent Theorem 24.}  Let $\pi=(d_1,d_2,\cdots,d_n) \in
  GS_n$, $n \geq 6,$ and $\pi\neq (2^{16})$, $(2^{17})$, $(4,
  2^{14})$, $(4,  2^{15})$, $(6,  2^{13})$, $(6,  2^{14})$,
  Then $\pi$ is potentially
$\sb 4 C\sb 6$-graphic if and only if $\pi$ be a potentially
$C_6$-graphic sequence, and  $\pi$ satisfies one of the following
conditions:
\par
  $(1)$ $d_3\geq 3$, and if $d_1 =d_3 =3, d_4 =2,$ then $d_{10} =2$;
  \par
  $(2)$ $d_2 \geq 4, d_3 =2, d_7 = 2$;
  \par
  $(3)$ $d_2 =3, d_3 =2,$ and if  $4 \geq d_1 \geq 3,$ then $d_{10} =2,$ and if  $n-4 \geq d_1 \geq 5,$ then $d_{9} =2$.
\par
  $(4)$ $d_2 =2,$ and if  $3 \geq d_1 \geq 2,$ then $d_{15} =2,$ and if  $5 \geq d_1 \geq 4,$ then $d_{14} =2$,
   and if  $n-7 \geq d_1 \geq 6,$ then $d_{13} =2$.
\par
 \textbf{\noindent Theorem 25.}  Let $\pi=(d_1,d_2,\cdots,d_n) \in
  GS_n$, $n \geq 6,$ and $\pi\neq (2^{12})$, $(2^{13})$, $(4,
  2^{10})$, $(4,  2^{11})$,
  Then $\pi$ is potentially
$\sb 5 C\sb 6$-graphic if and only if $\pi$ be a potentially
$C_6$-graphic sequence, and  $\pi$ satisfies one of the following
conditions:
\par
  $(1)$ $d_2\geq 3$;
  \par
  $(2)$ $3 \geq d_1 \geq 2, d_2 =2, d_11 = 2$;
  \par
  $(3)$   $n-6 \geq d_1 \geq 4, d_2 =2, d_{10} =2.$

\par
  Luo and Warner [54] characterized the potentially
 $K_4$-graphic sequences.
 \par
 \textbf{\noindent Theorem 26.}  Let $\pi=(d_1,d_2,\cdots,d_n)$
 be a graphic sequence without zero terms and with $d_4 \geq 3$ and $n\geq 4$. Then $\pi$ is potentially
$K_4$-graphic if and only if $d_4 \geq 3 $ and $\pi\neq (n-1, 3^s,
!^{n-s-1})$ for each $s =4,5$ except the following sequences:
\par
$n =5 :$ $(4, 3^4)$, $(3^4, 2)$;
\par
$n =6 :$ $(4^6)$, $(4^2, 3^4)$, $(4, 3^4, 2)$, $(3^6)$, $(3^5,1)$,
$(3^4, 2^2)$;
\par
$n =7 :$ $(4^7)$, $(4^3, 3^4)$, $(4, 3^6)$, $(4, 3^5, 1)$, $(3^6,
2)$, $3^5, 2, 1)$;
\par
$n =8 :$ $(3^7, 1)$, $(3^6, 1^2)$.
\par
   Eschen and Niu [14] and Lai[31] independently characterized the potentially
 $K_4-e$-graphic sequences.
 \par
 \par
\textbf{\noindent Theorem 27.}  Let $\pi=(d_1,d_2,\cdots,d_n)$
 be a graphic sequence with $n\geq 4$. Then $\pi$ is potentially
$(K_4-e)$-graphic if and only if the following conditions hold:
\par
  $(1)$ $d_2\geq 3$.
  \par
  $(2)$ $d_4 \geq 2$
\par
  $(3)$ If $n = 5, 6,$ then $\pi\neq (3^2, 2^{n-2})$ and $\pi\neq (3^6)$.
\par

 Yin and Yin [73] characterized the
 potentially $K_5-e$ and $K_6$-graphic sequences.
 \par
\textbf{\noindent Theorem 28}  Let $n\geq5$ and
$\pi=(d_1,d_2,\cdots,d_n)\in NS_n$
 be a positive graphic sequence with $d_3\geq4$ and $d_5\geq3$. Then $\pi$ is potentially
$K_5-e$-graphic if and only if $\pi$ is not one of the following
sequences: $(n-1,4^6,1^{n-7})$, $(n-1,4^2,3^4,1^{n-7})$,
$(n-1,4^2,3^3,1^{n-6})$;
\par
$n=6$: $(4^6), (4^4,3^2), (4^3,3^2,2)$;
\par
$n=7$: $(4^3,3^4)$, $(5^2,4,3^4)$, $(4^7)$, $(4^5,3^2)$,
$(5,4^3,3^3)$, $(5^2,4^5)$, $(5,4^5,3)$, $(4^3,3^2,2^2)$,
$(4^4,3^2,2)$, $(5,4^2,3^3,2)$, $(4^6,2)$, $(4^3,3^3,1)$;
\par
$n=8$: $(5^8)$, $(4^8)$, $(5^2,4^6)$, $(6,4^7)$, $(4^4,3^4)$,
$(5,4^2,3^5)$, $(4^6,3^2)$, $(5,4^6,3)$, $(4^3,3^4,2)$, $(4^7,2)$,
$(4^4,3^3,1)$, $(5,4^2,3^4,1)$, $(4^3,3^3,2,1)$, $(4^6,3,1)$,
$(5,4^6,1)$;
\par
$n=9$: $(4^9)$, $(4^3,3^5,1)$, $(4^8,2)$, $(4^7,3,1)$, $(5,4^7,1)$,
$(4^3,3^4,1^2)$, $(4^7,1^2)$;
\par
$n=10$: $(4^8,1^2)$.
\par
\textbf{\noindent Theorem 29}  Let $n\geq18$ and
$\pi=(d_1,d_2,\cdots,d_n)\in NS_n$
 be a positive graphic sequence with $d_6\geq5$. Then $\pi$ is potentially
$A_6$-graphic if and only if $\pi_6 \not \in$$ \{(2), (2^2), (3,1),
(3^2), (3,2,1), (3^2,2), (3^3,1), (3^2,1^2)\}$.
\par
  Yin and Chen [61] characterized the
 potentially $K_{r,s}$-graphic sequences for $r=2,s=3$ and
 $r=2,s=4$.\par
\textbf{\noindent Theorem 30} Let $n\geq5$ and
$\pi=(d_1,d_2,\cdots,d_n)\in GS_n$. Then $\pi$ is potentially
$K_{2,3}$-graphic if and only if $\pi$ satisfies the following
conditions:
\par
  $(1)$ $d_2\geq3$ and $d_5\geq2$;
  \par
  (2)If $d_1=n-1$ and $d_2=3$, then $d_5=3$;
  \par
 (3) $\pi\neq(3^2,2^4)$, $(3^2,2^5)$, $(4^3,2^3)$, $(n-1,3^5,1^{n-6})$
  and $(n-1,3^6,1^{n-7})$.
  \par
\textbf{\noindent Theorem 31} Let $n\geq6$ and
$\pi=(d_1,d_2,\cdots,d_n)\in GS_n$. Then $\pi$ is potentially
$K_{2,4}$-graphic if and only if $\pi$ satisfies the following
conditions:
\par
  $(1)$ $d_2\geq4$ and $d_6\geq2$;
  \par
  (2)If $d_1=n-1$ and $d_2=4$, then $d_3=4$ and $d_6\geq3$;
  \par
 (3) $\pi\neq(4^3,2^4)$, $(4^2,2^5)$, $(4^2,2^6)$, $(5^2,4,2^4)$, $(5^3,3,2^3)$,
 $(6,5^2,2^5)$, $(5^3,2^4,1)$, $(6^3,2^6)$, $(n-1,4^2,
 3^4,1^{n-7})$, $(n-1,4^2,3^5,1^{n-8})$, $(n-2,4^2,2^3,1^{n-6})$,
  and $(n-2,4^3,2^2,1^{n-6})$.
 \par
 Chen [3] characterized the
 potentially $K_5-2K_2$-graphic sequences for $5\leq n\leq8$.
 Hu and Lai [23] characterized the potentially
  $K_5-P_3$, $K_5-A_3$,
$K_5-K_3$, $K_5-K_{1,3}$ and $K_5-2K_2$-graphic
 sequences where $A_3$ is $P_2\cup K_2$.
 \par
\textbf{\noindent Theorem 32}  Let $\pi=(d_1,d_2,\cdots,d_n)$
 be a graphic sequence with $n\geq5$. Then $\pi$ is potentially
$K_5-P_3$-graphic if and only if the following conditions hold:
\par
  $(1)$ $d_1\geq4$, $d_3\geq3$ and $d_5\geq2$.
\par
  $(2)$ $\pi\neq (4,3^2,2^3)$, $(4,3^2,2^4)$ and $(4,3^6)$.
\par
\textbf{\noindent Theorem 33}  Let $\pi=(d_1,d_2,\cdots,d_n)$ be a
graphic sequence with $n\geq5$. Then $\pi$ is potentially
$K_5-A_3$-graphic if and only if the following conditions hold:
\par
  $(1)$ $d_4\geq3$ and $d_5\geq2$.
\par
  $(2)$ $\pi\neq(n-1,3^3,2^{n-k},1^{k-4})$ where $n\geq6$ and
$k=4,5,\cdots,n-2$, $n$ and $k$ have the same parity.
\par
  $(3)$ $\pi\neq (3^4,2^2),(3^6),(3^4,2^3),(3^6,2),(4,3^6),(3^7,1),(3^8),(n-1,3^5,1^{n-6})$
and $(n-1,3^6,1^{n-7})$.

\par
\textbf{\noindent Theorem 34}  Let $\pi=(d_1,d_2,\cdots,d_n)$
 be a graphic sequence with $n\geq5$. Then $\pi$ is potentially
$K_5-K_3$-graphic if and only if the following conditions hold:
\par
  $(1)$ $d_2\geq4$ and $d_5\geq2$.
\par
  $(2)$ $\pi\neq (4^2,2^4)$, $(4^2,2^5)$, $(4^3,2^3)$ and $(4^6)$.
\par
\textbf{\noindent Theorem 35}  Let $\pi=(d_1,d_2,\cdots,d_n)$ be a
graphic sequence with $n\geq5$. Then $\pi$ is potentially
$K_5-K_{1,3}$-graphic if and only if the following conditions hold:
\par
  $(1)$ $d_1\geq4$ and $d_4\geq3$.
\par
  $(2)$ $\pi\neq (4,3^4,2)$, $(4^6)$, $(4^2,3^4)$, $(4,3^6)$, $(4^7)$, $(4,3^5,1)$, $(n-1,3^4,1^{n-5})$
and $(n-1,3^5,1^{n-6})$.

\par
\textbf{\noindent Theorem 36}  Let $\pi=(d_1,d_2,\cdots,d_n)$ be a
graphic sequence with $n\geq5$. Then $\pi$ is potentially
$K_5-2K_2$-graphic if and only if the following conditions hold:
\par
  $(1)$ $d_1\geq4$ and $d_5\geq3$;
\par
  $(2)$ $$ \pi \neq \left\{
    \begin{array}{ll}(n-i,n-j,3^{n-i-j-2k}, 2^{2k},1^{i+j-2})\\  \mbox{ $n-i-j$ is even;}\\
    (n-i,n-j,3^{n-i-j-2k-1}, 2^{2k+1},1^{i+j-2})
  \\     \mbox{$n-i-j$ is odd.} \end{array} \right. $$ \ \ \ \ \ where $1\leq j\leq
     n-5$ and $0\leq k\leq [{{n-j-i-4}\over 2}]$.
\par
  $(3)$ $\pi\neq (4^2,3^4)$,  $(4,3^4,2)$, $(5,4,3^5)$, $(5,3^5,2)$, $(4^7)$,
        $(4^3,3^4)$, $(4^2,3^4,2)$,\par  $(4,3^6)$, $(4,3^5,1)$,$(4,3^4,2^2)$, $(5,3^7)$,
        $(5,3^6,1)$, $(4^8)$, $(4^2,3^6)$, $(4^2,3^5,1)$,\par $(4,3^6,2)$,
        $(4,3^5,2,1)$, $(4,3^7,1)$, $(4,3^6,1^2)$,
        $(n-1,3^5,1^{n-6})$ and \par $(n-1,3^6,1^{n-7})$.
\par
Hu and Lai [21] characterized the potentially
 $K_5-C_4$-graphic
 sequences.
 \par
  \textbf{Theorem 37}  Let $\pi=(d_1,d_2,\cdots,d_n)$ be a graphic sequence with
 $n\geq5$. Then $\pi$ is potentially $(K_5-C_4)$-graphic if and
 only if the following conditions hold:

  \par
 $(1)$ $d_1\geq4$.
  \par
  $(2)$ $d_5\geq2$.

  \par
   $(3)$ $\pi\neq((n-2)^2,2^{n-2})$ for $n\geq6$, where the symbol
   $x^y$ stands for $y$ consecutive terms $x$.

  \par
   $(4)$ $\pi\neq(n-k,k+i,2^i,1^{n-i-2})$ where
   $i=3,4,\cdots,n-2k$ and $k=1,2,\cdots,[\frac{n-1}{2}]-1$.

  \par
   $(5)$ If $n=6$, then $\pi\neq (4,2^5)$.

  \par
   $(6)$ If $n=7$, then $\pi\neq(4,2^6)$.
    \par
    Hu and Lai [22] characterized the potentially
 $K_5-Z_4$-graphic
 sequences.
 \par
  \textbf{Theorem 38}  Let $\pi=(d_1,d_2,\cdots,d_n)$ be a graphic sequence with
 $n\geq5$. Then $\pi$ is potentially $(K_5-Z_4)$-graphic if and
 only if the following conditions hold:

  \par
 $(1)$ $d_1\geq 4$, $d_2\geq 3$ and $d_4\geq 2$.
  \par
Hu, Lai and
 Wang[25] characterized the potentially $K_5-P_4$ and $K_5-Y_4$-graphic sequences  where $Y_4$
is a tree on 5 vertices and 3 leaves.
 \par
 \textbf{Theorem 39} Let $\pi=(d_1,d_2,\cdots,d_n)$
be a graphic sequence with $n\geq5$. Then $\pi$ is potentially
$K_5-P_4$-graphic if and only if the following conditions hold:
\par
(1) $d_2\geq3$.
\par
(2) $d_5\geq2$.
\par
(3) $\pi\neq(n-1,k,2^t,1^{n-2-t})$ where $n\geq5$,
$k,t=3,4,\cdots,n-2$, and, $k$ and $t$ have different parities.
\par
(4) For $n\geq5$, $\pi\neq(n-k,k+i,2^i,1^{n-i-2})$ where
$i=3,4,\cdots,n-2k$ and $k=1,2,\cdots,[\frac{n-1}{2}]-1$.
\par
(5) If $n=6,7$, then $\pi\neq (3^2,2^{n-2})$.
 \par
 \textbf{Theorem 40}
 Let $\pi=(d_1,d_2,\cdots,d_n)$ be a graphic sequence with $n\geq5$.
Then $\pi$ is potentially $K_5-Y_4$-graphic if and only if the
following conditions hold:

\par
(1) $d_3\geq3$.

\par
(2) $d_4\geq2$.

\par
(3) $\pi\neq(3^6)$.
 \par
   Hu and Lai [24] characterized the potentially
 $K_{3,3}$ and $K_6-C_6$-graphic
 sequences.
 \par
\textbf{\noindent Theorem 41}  Let $\pi=(d_1,d_2,\cdots,d_n)$
 be a graphic sequence with $n\geq6$. Then $\pi$ is potentially
$K_{3,3}$-graphic if and only if the following conditions hold:
\par
  $(1)$ $d_6\geq3$;\par
  $(2)$ For $i=1,2$, $d_1=n-i$ implies $d_{4-i}\geq4$;\par
  $(3)$ $d_2=n-1$ implies $d_3\geq5$ or $d_6\geq4$;\par
  $(4)$ $d_1+d_2=2n-i$ and $d_{n-i+3}=1 (3\leq i\leq n-4)$ implies $d_3\geq5$ or $d_6\geq4$;\par
  $(5)$ $d_1+d_2=2n-i$ and $d_{n-i+4}=1 (4\leq i\leq n-3)$ implies $d_3\geq4$;\par
  $(6)$ $\pi=(d_1,d_2,3^4,2^t,1^{n-6-t})$ or $(d_1,d_2,4^2,3^2,2^t,1^{n-6-t})$ implies $d_1+d_2\leq n+t+2$;\par
  $(7)$ $\pi=(d_1,d_2,4,3^4,2^t,1^{n-7-t})$ implies $d_1+d_2\leq n+t+3$;\par
  $(8)$ For $t=5,6$, $\pi\neq(n-i,k+i,4^t,2^{k-t},1^{n-2-k})$ where
  $i=1,\cdots,[{{n-k}\over2}]$ and $k=t,\cdots,n-2i$;\par
  $(9)$ $\pi\neq(5^4,3^2,2)$, $(4^6)$, $(3^6,2)$, $(6^4,3^4)$,
  $(4^2,3^6)$, $(4,3^6,2)$, $(3^6,2^2)$, $(3^8)$, $(3^7,1)$,
  $(4,3^8)$, $(4,3^7,1)$, $(3^8,2)$, $(3^7,2,1)$, $(3^9,1)$,
  $(3^8,1^2)$, $(n-1,4^2,3^4,1^{n-7})$, $(n-1,4^2,3^5,1^{n-8})$,
  $(n-1,5^3,3^3,1^{n-7})$,
  $(n-2,4,3^5,1^{n-7})$, $(n-2,4,3^6,1^{n-8})$, $(n-3,3^6,1^{n-7})$, $(n-3,3^7,1^{n-8})$.
\par
\par
\textbf{\noindent Theorem 42} Let $\pi=(d_1,d_2,\cdots,d_n)$
 be a graphic sequence with $n\geq6$. Then $\pi$ is potentially
$K_6-C_6$-graphic if and only if the following conditions hold:
\par
  $(1)$ $d_6\geq3$;\par
  $(2)$ For $i=1,2$, $d_1=n-i$ implies $d_{4-i}\geq4$;\par
  $(3)$ $d_2=n-1$ implies $d_4\geq4$;\par
  $(4)$ $d_1+d_2=2n-i$ and $d_{n-i+3}=1 (3\leq i\leq n-4)$ implies $d_4\geq4$;\par
  $(5)$ $d_1+d_2=2n-i$ and $d_{n-i+4}=1 (4\leq i\leq n-3)$ implies $d_3\geq4$;\par
  $(6)$ $\pi=(d_1,d_2,d_3,3^k,2^t,1^{n-3-k-t})$ implies $d_1+d_2+d_3\leq n+2k+t+1$;\par
  $(7)$ $\pi=(d_1,d_2,3^4,2^t,1^{n-6-t})$ implies $d_1+d_2\leq n+t+2$;\par
  $(8)$ $\pi\neq(n-i,k,t,3^t,2^{k-i-t-1},1^{n-2-k+i})$ where $i=1,\cdots,[{{n-t-1}\over
  2}]$ and $k=i+t+1,\cdots,n-i$ and $t=4,5,\cdots,k-i-1$;\par
  $(9)$ $\pi\neq(3^6,2)$, $(4^2,3^6)$, $(4,3^6,2)$, $(3^6,2^2)$, $(3^8)$, $(3^7,1)$,
  $(4,3^8)$, $(4,3^7,1)$, $(3^8,2)$, $(3^7,2,1)$, $(3^9,1)$,
  $(3^8,1^2)$, $(n-1,4^2,3^4,1^{n-7})$, $(n-1,4^2,3^5,1^{n-8})$,
  $(n-2,4,3^5,1^{n-7})$, $(n-2,4,3^6,1^{n-8})$, $(n-3,3^6,1^{n-7})$, $(n-3,3^7,1^{n-8})$.
\par
Xu and Hu[57] characterized the potentially $K_{1,4}+e$-graphic
sequences. Chen and Li [8] characterized the potentially
$K_{1,t}+e$-graphic sequences.
\par
\textbf{\noindent Theorem 43} Let $\pi=(d_1,d_2,\cdots,d_n)$
 be a graphic sequence with $n\geq5$. Then $\pi$ is potentially
$K_{1,4}+e$-graphic if and only if $d_1\geq4, d_3\geq2.$
\par
\textbf{\noindent Theorem 44} Let $t\geq3$,
$\pi=(d_1,d_2,\cdots,d_n)$
 is a graphic sequence with $n\geq t+1$. Then $\pi$ is potentially
$K_{1,t}+e$-graphic if and only if $d_1\geq t, d_3\geq2.$
\par
\section*{Open Problems}
\par
\par
Problem 1. Determine $\sigma(K_{r+1}-G, n)$ for $G$ is a graph on
$k$ vertices and $j$ edges which
    contains a graph $K_{3} \bigcup K_{1,3}$ but
     not contains a cycle on $4$ vertices and not contains $Z_4$, $P_3$.
\par
Problem 2. Determine $\sigma(K_{r+1}-G, n)$ for $G= K_{3} \bigcup
iK_{2}\bigcup jP_{2}\bigcup tK_{3}$.

\par
Problem 3. Determine $\sigma(K_{r+1}-G, n)$ for graph $G$ which
contains $C_3$, $C_4,\cdots, C_l$ but not contains a cycle on $l+1$
vertices($4\leq l\leq r$).
\par
 Problem 4. Determine $\sigma(K_{r+1}-G, n)$
for graph $G$ which contains $C_3$, $C_4,\cdots, C_{r+1}$.
\par

 Problem 5.
Determine $\sigma(K_{r+1}-G, n)$ for small $n$.
\par
 Problem 6.
Characterize potentially $K_{r+1}-G$-graphic sequences for the
remaining $G$.
\par

\par
 \section*{Acknowledgment}
 The first author is particularly indebted to Professor Jiongsheng Li for
 introducing him to degree sequences. The authors wish to thank Professor Gang Chen, R.J.
Gould,  Jiongsheng Li, Rong Luo,  John R. Schmitt,  Zi-Xia Song,
Amitabha Tripathi, Jianhua Yin and Mengxiao Yin for sending some
their papers to us.

\par

\end{document}